\documentclass[12pt]{amsart}
\usepackage{amsmath}
\usepackage{amssymb}
\usepackage{bm}
\usepackage{graphicx}
\usepackage{verbatim}
\newtheorem{Theorem}{Theorem}
\newtheorem{Lemma}{Lemma}

\begin{document}

\title[Localized Quasimodes in Large Subspaces]{Partially Localized Quasimodes in Large Subspaces\\
(PRELIMINARY VERSION -- COMMENTS WELCOME)}
\author{Shimon Brooks}
\thanks{The author was partially supported by NSF grant DMS-1101596.}

\maketitle

{\em Abstract:}  We consider spaces of high-energy quasimodes for the Laplacian on a compact hyperbolic surface, and show that when the spaces are large enough, one can find quasimodes that exhibit strong localization phenomena.  Namely, take any constant $c$, and a sequence of $cr_j$-dimensional spaces $\mathcal{S}_j$ of quasimodes, where $\frac{1}{4}+r_j^2\to\infty$ is an approximate eigenvalue for $\mathcal{S}_j$.  Then we can find a sequence of vectors $\psi_j\in\mathcal{S}_j$, such that any weak-* limit point of the microlocal lifts of $|\psi_j|^2$ localizes a positive proportion of its mass on a singular set of codimension $1$.  This result is sharp, in light of the QUE result of \cite{jointQmodes} for certain joint quasimodes that include spaces of size $o(r_j)$, with arbitrarily slow decay.
\medskip

\section{Introduction}

The Quantum Unique Ergodicity (QUE) Conjecture of Rudnick-Sarnak \cite{RS} states that  eigenfunctions of the Laplacian on Riemannian manifolds of negative sectional curvature become 
equidistributed in the high-energy limit.  Although there exist so-called ``toy models" of quantum chaos that do not exhibit this behavior (see eg. \cite{FNDB, ANb, KelPert}), it has been suggested that large degeneracies of the quantum propagator may be responsible for some of
these phenomena
(see eg. \cite{SarnakProgress}).  Since the Laplacian on a surface of negative curvature is not expected to have large degeneracies, one can explore this aspect and introduce ``degeneracies" by considering quasimodes, or approximate eigenfunctions, in place of true eigenfunctions--- relaxing the order of approximation to true eigenfunctions yields larger spaces of quasimodes, mimicking higher-dimensional eigenspaces.  Studying the properties of such quasimodes--- and, especially, the effect on equidistribution of varying the order of approximation--- can help shed light on the overall role of spectral degeneracies in the theory.
 
 Let $X=\Gamma\backslash SL(2,\mathbb{R})=S^*M$ compact.  We normalize the uniform measure $dz$ on $M$ (and the measure $dx$ on $X$) to have total volume $1$.
We define an $\omega(r)$-{\bf quasimode with approximate parameter} $r$ to be a function $\psi$ satisfying
$$||(\Delta + (\frac{1}{4}+r^2))\psi||_2 \leq r \omega(r)||\psi||_2$$
The  factor of $r$ in our definition comes from the fact that $r$ is essentially the square-root of the Laplace eigenvalue.  
For any constant $C$, denote by $S_C(r)$ the space spanned by eigenfunctions of spectral parameter in $[r-C, r+C]$; then for large $r$ (in particular, $r>C$), any vector in $S_C(r)$ is a $3C$-quasimode of approximate parameter $r$.

For a sequence of spectral parameters $r_j\to\infty$, consider the space $S_C(r_j)\subset C^\infty(\Gamma\backslash\mathbb{H})$, spanned by the eigenfunctions whose spectral parameters lie in $[r_j-C, r_j+C]$; coarse estimates on the remainder in Weyl's Law show that this space has dimension $\dim{S_C(r_j)} \gtrsim C r_j$.  To any vector $\psi\in S_C(r_j)$ we associate a measure $|\Psi|^2dVol$ on $S^*M$ (see section~\ref{ml lifts construction}) called the {\bf microlocal lift} of $|\psi|^2dz$.  A sequence $\{\psi_j\in S_C(r_j)\}$ is said to satisfy the QUE property if the measures $|\Psi_j|^2dVol$ converge weak-* to the uniform (Liouville) measure on $S^*M$.

In joint work with Lindenstrauss \cite{jointQmodes}, we studied certain cases of {\em  joint quasimodes}--- eg., of the Laplacian and one Hecke operator--- and found that QUE held for all sequences of functions, that were jointly $o(1)$-quasimodes for both the Laplacian and one Hecke operator.  It is important to note that there are spaces of such joint $o(1)$ quasimodes of size $o(r)$ with arbitrarily slow decay; these are considerably larger than the spaces of Laplace-quasimodes that are expected to satisfy QUE without any Hecke assumption.  This is a testament to the rigidity imposed by the additional structure of the Hecke correspondence, as already apparent in \cite{Lin}.

Here we show that, in fact, these $o(r)$-dimensional spaces are as large as possible for QUE.

\begin{Theorem}\label{optimality}
Let $\mathcal{S}_j\subset S_C(r_j)$ be a subspace of dimension $c r_j$ for each $S_C(r_j)$.  Then there exists a sequence of quasimodes $\psi_j\in \mathcal{S}_j$, such that any weak-* limit point of the microlocal lifts $|\Psi_j|^2dVol$ gives positive measure to a subset of $X$ of codimension $1$.
\end{Theorem}
The results of \cite{jointQmodes} show that for these special arithmetic manifolds, there can be subspaces of $S_C(r_j)$, of dimension $o(r_j)$, such that any sequence of quasimodes taken from these spaces must satisfy QUE.  Theorem~\ref{optimality} shows that if the subspaces of quasimodes are taken to be any larger--- i.e., of dimension $\geq cr_j$ for some fixed constant $c$--- then one can always find bad sequences that do not satisfy QUE.  In this sense, the special joint $o(1)$-quasimodes of \cite{jointQmodes} are optimally degenerate for QUE.

The codimension $1$ subset of $X=\Gamma\backslash SL(2,\mathbb{R})$ we have in mind is the collection of geodesic segments through a single base point $p\in\Gamma\backslash\mathbb{H}$; i.e., cotangent vectors pointing radially towards or away from $p$, at distance $\leq \tau$, where $\tau$ is a fixed number depending only on the manifold $M$.  For each $j$, we will find a base point $p_j\in\Gamma\backslash\mathbb{H}$ and a quasimode $\psi_j\in \mathcal{S}_j$ that is large at $p_j$, and thus the microlocal lift $|\Psi_j|^2dVol$ (see section~\ref{ml lifts}) will be enhanced on vectors pointing radially relative to $p_j$.  We will then take $p$ to be a limit point of $\{p_j\}$ in the compact manifold $\Gamma\backslash\mathbb{H}$.

The idea is to use a kernel that is spectrally localized near $S_C(r_j)$, and spatially localized near radial vectors around $p_j$; we describe such a kernel in section~\ref{ml kernel}.  Since the kernel is spectrally localized near $S_C(r_j)$, it will strongly correlate with our $\Psi_j$ (up to factors depending on $c$, $C$, and the manifold $M$), and show that the latter must also localize a fixed positive proportion of its mass near our codimension $1$ subset.

{\bf Acknowledgements.}  This paper was motivated by joint work \cite{jointQmodes} with Elon Lindenstrauss, whom we thank for encouragement and many helpful discussions.

\section{Microlocal Lifts of Quasimodes}\label{ml lifts}

\subsection{Some Harmonic Analysis on $SL(2,\mathbb{R})$}\label{harmonic analysis}
We begin by reviewing some harmonic analysis on $SL(2,\mathbb{R})$ that we will need.  
Throughout, we write $X=\Gamma\backslash SL(2,\mathbb{R})$ and $M=\Gamma\backslash\mathbb{H} = \Gamma\backslash SL(2,\mathbb{R})/K$, where $K=SO(2)$ is the maximal compact subgroup.

Fix an orthonormal basis $\{\phi_l\}$ of $L^2(M)$ consisting of Laplace eigenfunctions, which we can take to be real-valued for simplicity.  Each eigenfunction generates, under right translations, an irreducible representation $V_{l} = \overline{\{\phi_l(xg^{-1}) : g\in SL(2,\mathbb{R})\}}$ of $SL(2,\mathbb{R})$, which span a dense subspace of $L^2(X)$.  

We distinguish the pairwise orthogonal {\bf weight spaces} $A_{2n}$  in each representation, consisting of those functions satisfying $f(xk_\theta) = e^{ 2i n\theta}f(x)$ for all $k_\theta = \begin{pmatrix} \cos{\theta} & \sin{\theta}\\ -\sin{\theta} & \cos{\theta}\end{pmatrix} \in K$ and $x\in X$.  The weight spaces together span a dense subspace of $V_{l}$.  Each weight space is one-dimensional in $V_{l}$, spanned by $\phi^{(l)}_{2n}$ where 
\begin{eqnarray*}
\phi_0^{(l)} & = & \phi_l \in A_0\\
(ir_l + \frac{1}{2}+n)\phi_{2n+2}^{(l)} & = & E^+\phi_{2n}^{(l)}\\ 
(ir_l + \frac{1}{2} -n)\phi_{2n-2}^{(l)} & = & E^-\phi_{2n}^{(l)}
\end{eqnarray*}
Here $E^+$ and $E^-$ are the {\bf raising} and {\bf lowering operators}, first-order differential operators corresponding to $\begin{pmatrix} 1 & i\\ i & -1\end{pmatrix}\in \mathfrak{sl}(2,\mathbb{R})$ and $\begin{pmatrix} 1 & -i\\ -i & -1 \end{pmatrix}\in \mathfrak{sl}(2,\mathbb{R})$ in the Lie algebra.  
The normalized pseudodifferential operators 
\begin{eqnarray*}
R & = & \frac{E^+}{ir_l+\frac{1}{2}+n} : \phi_{2n}^{(l)}\mapsto \phi_{2n+2}^{(l)}\\
R^{-1} & = & \frac{E^-}{ir_l+\frac{1}{2}-n} : \phi_{2n}^{(l)}\mapsto \phi_{2n-2}^{(l)}
\end{eqnarray*}
are unitary and left-invariant, and each $\phi_{2n}$ is a unit vector.
We define the distribution
$$\Phi_\infty^{(l)} = \sum_{n=-\infty}^{\infty} \phi_{2n}^{(l)}$$
and extend this definition by linearity to  $\Psi_\infty = \sum_{n=-\infty}^\infty \psi_{2n} = \sum_{n=-\infty}^\infty R^n\psi$ for $\psi\in S_C(r)$, where each $\psi_{2n} = R^n\psi$.

\vspace{.2in}
Let $k$ be a $K$-bi-invariant function on $SL(2,\mathbb{R})$--- i.e., a radial function on $\mathbb{H}$.  Then any Laplace eigenfunction $\phi$ of eigenvalue $\frac{1}{4}+s^2$ is also an eigenfunction of convolution with $k$, with eigenvalue given by the spherical transform $h(s)$ of $k$.  The spherical transform is related to $k$ by the Selberg/Harish-Chandra transform  (see eg. \cite[Chapter 1.8]{Iwaniec})
\begin{eqnarray}
h(s) & = & \int_{-\infty}^\infty e^{isu}g(u)du\nonumber\\
g(u) & = & 2Q\left(\sinh^2\left(\frac{u}{2}\right)\right)\nonumber\\
k(t) & = & -\frac{1}{\pi} \int_t^\infty \frac{dQ(\omega)}{\sqrt{\omega-t}}\label{Selberg/HC}
\end{eqnarray}
The coordinate $t(z,w) = 2\sinh^2(d(z,w)/2)$ is often more convenient for calculations ($dt$ is the radial volume measure on $\mathbb{H}$).  What is most important for our purposes is that whenever $g$ is compactly supported in the interval $[-\tau, \tau]$, the kernel $k$ will be supported in the ball of radius $\tau$ in $\mathbb{H}$.

We can write such a $k$ as a (left-$K$-invariant) function on $\mathbb{H}$, and use Helgason's Fourier inversion \cite{Hel} to write
$$k(z) = \int_{s=0}^\infty \int_{B} e^{(is+\frac{1}{2})<z,b>}\hat{k}(s,b) s\tanh{(\pi s)} dsdb$$
where $b\in B$ runs over the boundary $S^1$ of the disc model for $\mathbb{H}$, and $<z,b>$ represents the (signed) distance to the origin from the horocycle through the point $z\in\mathbb{H}$ tangent to $b\in B$.  Since each plane wave $e^{(-is+\frac{1}{2})<\cdot,b>}$ is an eigenfunction of spectral parameter $s$, the Fourier transform
\begin{eqnarray*}
\hat{k}(s,b) & = & \int_\mathbb{H} e^{(-is+\frac{1}{2})<z,b>}k(z)dz\\
& = & h(s)
\end{eqnarray*}
so that
$$k(z) = \int_{s=0}^\infty \left(  \int_{B} e^{(is+\frac{1}{2})<z,b>}db \right) h(s) s\tanh{(\pi s)} ds$$

It will be more convenient  to write as in \cite{Z2} 
$$e^{(is+\frac{1}{2})<z,b>}db = e^{(is-\frac{1}{2})<z,b>}d\theta =  e^{(is-\frac{1}{2})\varphi(g.k_\theta)}d\theta$$
where $\varphi(g)$ is the signed distance from the origin to the horocycle through $g\in SL(2,\mathbb{R})$, and $k_\theta$ parametrizes the $SO(2)$ fibre  $gK$.  Since $\varphi$ is left $K$-invariant and right $N$-invariant, it is convenient to use $KAN$ coordinates to write
$$g=\begin{pmatrix}  \cos{\theta} & -\sin{\theta}\\ \sin{\theta} & \cos{\theta}\end{pmatrix} \begin{pmatrix} e^{t/2} & 0\\ 0 & e^{-t/2}\end{pmatrix}  \begin{pmatrix}  1 & n\\ 0 & 1\end{pmatrix} = \begin{pmatrix} a& b\\ c& d\end{pmatrix}$$
so that this distance is given by $\varphi(g) = t =\log{(a^2+c^2)}$.

\subsection{Construction of the Microlocal Lifts}\label{ml lifts construction}
We set 
$$I_\psi(f) =\langle Op(f)\psi, \psi\rangle := \langle f \Psi_\infty, \psi\rangle = \lim_{N\to\infty} \left\langle f \sum_{n=-N}^N \psi_{2n}, \psi_0\right\rangle$$
according to the pseudo-differential calculus of \cite{Z2}, which clearly restricts to the measure $|\psi(z)|^2dz$ when applied to $K$-invariant functions $f\in C^\infty(M)$, by orthogonality of the weight spaces.
Note moreover that this limit is purely formal for $K$-finite $f$, and since these $K$-finite functions are dense in the space of smooth functions, we can restrict our attention to these.  We denote by $A_{2n}$ the $n$-th weight space, consisting of smooth functions that transform via $f(xk_\theta) = e^{2i n\theta} f(x)$ for all $x\in X$.

\begin{Lemma}\label{asymp to op}
Let $\psi\in S_C(r)$ be a unit vector, and set 
$$\Psi := \sqrt{\frac{3L}{2L^2+1}}\sum_{|n|\leq\sqrt{r}} \frac{L - |n|}{L} \psi_{2n}$$  
for $L:= \lfloor \sqrt{r}\rfloor$.  
Then for any $K$-finite $f\in \sum_{n=-N_0}^{N_0}A_{2n}$, we have
$$I_\psi(f) = \langle f\Psi, \Psi \rangle +O_{f,C}(r^{-1/2})$$
\end{Lemma}

Note that the prefactor $\sqrt{\frac{3L}{2L^2+1}}\sim \sqrt{\frac{3}{2}}r^{-1/4}$ is simply an $L^2$-normalization of the Fej\'er coefficients $\frac{L - |n|}{L} $.
The  proof of the Lemma is  identical to that of  \cite{jointQmodes}, the only difference being this use of  Fej\'er coefficients in place of the Dirichlet coefficients $\frac{1}{\sqrt{2L+1}}$ used in \cite{LinHxH} and \cite{jointQmodes}.  The extra smoothness provided by using Fej\'er coefficients  will be exploited in section~\ref{ml kernel}.

{\em Proof:}  First, we wish to show that
$$\langle f \psi_{2n}, \psi_{2m}\rangle = \langle f \psi_{2n+2}, \psi_{2m+2}\rangle +O_{f,C}(r^{-1})$$
for all $-\sqrt{r}\leq n,m \leq \sqrt{r}$, satisfying $|n-m|\leq N_0$ (if the latter condition is not met, both inner products are trivial, by orthogonality of the weight spaces).  We will work individually with each pair of spectral components of $\psi$, and then re-average over the spectral decomposition; therefore, we write $\psi^{(r_1)}$ and $\psi^{(r_2)}$ for the projections of $\psi$ to the eigenspaces of parameters $r_1$ and $r_2$, respectively.  Recall that $r_1, r_2 = r+O_C(1)$ by the condition $\psi\in S_C(r_j)$.

We have 
\begin{eqnarray*}
\lefteqn{\langle f \psi^{(r_1)}_{2n}, \psi^{(r_2)}_{2m}\rangle }\\
& = & \frac{1}{(ir_1-n-\frac{1}{2})(-ir_2 -m-\frac{1}{2})}\langle f E^-\psi^{(r_1)}_{2n+2}, E^-\psi^{(r_2)}_{2m+2}\rangle\\
& = & \frac{1}{(ir_1-n-\frac{1}{2})(-ir_2 -m-\frac{1}{2})}\Big( \langle E^-(f\psi^{(r_1)}_{2n+2}), E^-\psi_{2m+2}^{(r_2)}\rangle - \langle E^-(f)\psi^{(r_1)}_{2n+2}, E^-\psi_{2m+2}^{(r_2)}\rangle \Big)\\
& = & - \frac{1}{(ir_1-n-\frac{1}{2})(-ir_2 -m-\frac{1}{2})} \langle f\psi^{(r_1)}_{2n+2}, E^+E^- \psi^{(r_2)}_{2m+2}\rangle + O(r_1^{-1}) \langle E^-(f)\psi^{(r_1)}_{2n+2}, \psi^{(r_2)}_{2m}   \rangle\\
& = & - \frac{(-ir_2+m+\frac{1}{2})(-ir_2-m-\frac{1}{2})}{(ir_1-n-\frac{1}{2})(-ir_2 -m-\frac{1}{2})}
\langle f\psi^{(r_1)}_{2n+2}, \psi^{(r_2)}_{2m+2}\rangle + O(r_1^{-1}) \langle E^-(f)\psi^{(r_1)}_{2n+2}, \psi^{(r_2)}_{2m}   \rangle\\
& = & \langle f\psi_{2n+2}^{(r_1)}, \psi_{2m+2}^{(r_2)}\rangle + \frac{i(r_2-r_1) + (n-m)}{ir_1-n-\frac{1}{2}} \langle f\psi_{2n+2}^{(r_1)}, \psi_{2m+2}^{(r_2)}\rangle + O(r_1^{-1}) \langle E^-(f)\psi^{(r_1)}_{2n+2}, \psi^{(r_2)}_{2m}   \rangle
\end{eqnarray*}

We now average over $r_2$, and since $i(r_2-r_1) + (n-m) = O_{f,C}(1)$ (recall that $|n-m|\leq N_0 = O_f(1)$), we have
\begin{eqnarray*}
\lefteqn{\langle f \psi_{2n}^{(r_1)}, \psi_{2m}\rangle - \langle f\psi_{2n+2}^{(r_1)}, \psi_{2m+2}\rangle}\\
& = & \frac{1}{ir_1-n-\frac{1}{2}} \left\langle f \psi_{2n+2}^{(r_1)}, \sum_{r_2} O_{f,C}(1) \psi_{2m+2}^{(r_2)}\right\rangle + O(r_1^{-1}) \left\langle E^-(f) \psi_{2n+2}^{(r_1)}, \psi_{2m}\right\rangle
\end{eqnarray*}
and further averaging over $r_1$ gives 
\begin{eqnarray*}
\lefteqn{\langle f\psi_{2n}, \psi_{2m}\rangle - \langle f\psi_{2n+2}, \psi_{2m+2}\rangle} \\
& = &  \left\langle f \sum_{r_1} O(r_1^{-1}) \psi_{2n+2}^{(r_1)}, \sum_{r_2} O_{f,C}(1) \psi_{2m+2}^{(r_2)}\right\rangle + \left\langle E^-(f) \sum_{r_1} O(r_1^{-1})  \psi_{2n+2}^{(r_1)}, \psi_{2m}\right\rangle
\end{eqnarray*}

Now the error terms are estimated by Cauchy-Schwarz, giving
\begin{eqnarray*}
\lefteqn{\left\langle f \sum_{r_1} O(r_1^{-1})  \psi_{2n+2}^{(r_1)}, \sum_{r_2} O_{f,C}(1) \psi_{2m+2}^{(r_2)}\right\rangle}\\ 
& \leq & ||f||_\infty \left|\left|\sum_{r_1}  O(r_1^{-1})\psi_{2n+2}^{(r_1)}\right|\right|_2 \left|\left|  \sum_{r_2} O_{f,C}(1) \psi_{2m+2}^{(r_2)}\right|\right|_2\\
& \leq & O_{f,C}(r_1^{-1}) ||\psi_{2n+2}||_2 ||\psi_{2m+2}||_2\\
& = & O_{f,C}( r_1^{-1})
\end{eqnarray*}
and
\begin{eqnarray*}
\langle E^-(f) \sum_{r_1} O(r_1^{-1}) \psi_{2n+2}^{(r_1)}, \psi_{2m}\rangle 
& \leq &
||E^-(f)||_\infty O(r_1^{-1})||\psi_{2n+2}||_2 ||\psi_{2m}||_2\\
& = & O_f(r_1^{-1})
\end{eqnarray*}
by using the orthogonality of the $\psi_{2n+2}^{(r_1)}$ to estimate
\begin{eqnarray*}
\left|\left|	\sum_{r_1} O(r_1^{-1})\psi_{2n+2}^{(r_1)}	\right|\right|_2^2 & = & O(r_1^{-2}) \sum_{r_1}||\psi_{2n+2}^{(r_1)}||_2^2\\
& = & O(r_1^{-2})||\psi_{2n+2}||_2^2
\end{eqnarray*}
and similarly
$$\left|\left|	\sum_{r_2} O_{f,C}(1)\psi_{2m+2}^{(r_2)}	\right|\right|_2^2
= O_{f,C}(1)||\psi_{2m+2}||_2^2$$

Therefore
\begin{eqnarray*}
\langle f\psi_{2n}, \psi_{2m}\rangle 
& = &  \langle f\psi_{2n+2}, \psi_{2m+2}\rangle + O_{f,C}(r^{-1})
\end{eqnarray*}
since  $|r_1- r| \leq C$ implies that $O_{f,C}(r_1^{-1})= O_{f,C}(r^{-1})$.  

We iterate this $|m|\leq \sqrt{r}$ times, arriving at
\begin{equation}\label{after iteration}
\langle f\psi_{2n}, \psi_{2m}\rangle = \langle f\psi_{2(n-m)}, \psi_0\rangle 
+ O_{f,C}(\sqrt{r}r^{-1})
\end{equation}
Now, by definition
\begin{eqnarray*}
\langle f\Psi, \Psi\rangle & = & \sum_{|m|,|n|\leq L} \frac{3(L - |n|)(L - |m|)}{2L^3 +L}\langle f\psi_{2n},\psi_{2m}\rangle\\
& = & \sum_{n=-L}^L \sum_{m=n-N_0}^{n+N_0} \frac{3(L - |n|)(L - |m|)}{2L^3 +L}\langle f\psi_{2n},\psi_{2m}\rangle
 - \sum_{|n|=L-N_0}^L O_f(N_0L^{-1})
 \end{eqnarray*}
 since each term $\frac{3(L-|n|)(L-|m|)}{2L^3 + L}\langle f \psi_{2n}, \psi_{2m}\rangle = O_f(L^{-1}))$, and for each value of $n$, there are at most $N_0$ values of $m$ such that the inner product is not trivial.  Thus, since $N_0=O_f(1)$, we get
 \begin{eqnarray*}
\langle f\Psi, \Psi\rangle & = &  \sum_{n=-L}^L \sum_{m=n-N_0}^{n+N_0} \frac{3 (L - |n|)(L - |m|)}{2L^3+L} \langle f\psi_{2n},\psi_{2m}\rangle + O_f(L^{-1})\\
& = & \frac{3}{2L^3 +L}\sum_{n=-L}^L \sum_{m=n-N_0}^{n+N_0} \Big((L - |n|)^2 - O(L|n-m|) \Big)\langle f\psi_{2n},\psi_{2m}\rangle + O_f(L^{-1})\\
& = & \left(\frac{3}{2L^3 +L}\sum_{n=- L}^L (L - |n|)^2\right) \sum_{n-m = -N_0}^{N_0}\Big(\langle f\psi_{2(n-m)},\psi_{0}\rangle +O_{f,C}(r^{-1/2})\Big) + O_f(L^{-1})
\end{eqnarray*}
by (\ref{after iteration}).  Therefore setting $l=n-m$  we finally obtain
\begin{eqnarray*}
\langle f\Psi, \Psi\rangle & = & \sum_{l=-N_0}^{N_0} \langle f\psi_{2l},\psi_{0}\rangle + O_{f,C}(r^{-1/2}) + O_f(L^{-1}) \\
& = & I_{\psi_j}  + O_{f,C}(r^{-1/2})
\end{eqnarray*} 
since $L\leq r^{-1/2}$.  $\Box$

For any given sequence $\{\psi_j\in S_C(r_j)\}$, we have constructed a sequence $\{\Psi_{j}\}_{j=1}^\infty$ such that the microlocal lifts $|\Psi_{j}|^2dVol$ are asymptotically equivalent to the distributions $I_{\psi_j}$.  It is these measures that we wish to study.

\section{A ``microlocal" kernel}\label{ml kernel}

Pick once and for all an orthonormal basis $\{\phi_l\}$ of $L^2(M)$ consisting of real-valued eigenfunctions.  Throughout, we will allow all implied constants to depend on $C$, $c$, and the manifold $M$.

We will need an auxilliary spherical kernel $k$.  Observe that the hypotheses of Theorem~\ref{optimality} are weaker  when $C$ is larger and $c$ smaller; in particular, we may assume that $C$ is sufficiently large.  We begin with
$$\tilde{h}(s) = r_j^{-1/2} \frac{\cosh{s/C}\cosh{r_j/C}}{\cosh{2s/C}+\cosh{2r_j/C}}$$
a scaled version of the spherical transform used in \cite{IwaniecSarnak}.  
The Fourier transform of $h$ is
$$\tilde{g}(\xi) =  \frac{C}{4} r_j^{-1/2} \frac{\cos(\xi r_j)}{\cosh\left(\frac{C\pi\xi}{2}\right)}$$

It will be convenient to cutoff $\tilde{g}$, so that our kernel will have compact support.  Let $\tau$ be a small--- but fixed--- number less than the radius of injectivity of $M$, and pick a smooth, non-negative, even cutoff function $\chi\in C_c^\infty([-\tau, \tau])$ whose Fourier transform $\hat{\chi}$ is also non-negative; we normalize so that $\hat{\chi}(s)\geq 1$ for all $|s|\leq 1$.  Note that $\chi$ depends only on $\tau$, and is independent of the parameters $C$ and $c$, as well as $r_j$.
Now let
$$g(\xi) = \tilde{g}(\xi)\chi(\xi) \in C_c^\infty([-\tau, \tau])$$
whereby the corresponding kernel $k$ given by (\ref{Selberg/HC}) will be compactly supported inside the ball of radius $\tau$ in $\mathbb{H}$.  

The spherical transform 
$h = \tilde{h} \ast \hat{\chi}$ satisfies
\begin{eqnarray}
||h||_\infty & \leq & ||\hat{\chi}||_{L^1} \cdot ||\tilde{h}||_\infty \nonumber\\
& \lesssim_\tau & r_j^{-1/2}\label{h small}
\end{eqnarray}
Moreover, since $\hat{\chi}\geq 1$ on $[-1,1]$, we have for $|s-r_j|\leq C$ 
\begin{eqnarray}
\min_{|s-r_j|\leq C} h(s) & \geq & \inf_{|s'-s|\leq 1} \tilde{h}(s) \geq \inf_{|s'-r_j|\leq (C+1)} r_j^{-1/2} \frac{\cosh{\frac{s'}{C}\cosh{\frac{r_j}{C}}}}{\cosh{\frac{2s'}{C}}+\cosh{\frac{2r_j}{C}}}\nonumber\\
& \geq &  \inf_{|s'-r_j|\leq (C+1)} r_j^{-1/2} \frac{1}{4} \frac{e^{s'/C+r_j/C}}{2\max\{e^{2s'/C}, e^{2r_j/C}\}}\nonumber\\
& \geq & \inf_{|s'-r_j|\leq (C+1)} \frac{1}{8}r_j^{-1/2} e^{-|s'-r_j|/C}\nonumber\\
& \geq & \frac{1}{100}r_j^{-1/2}\label{big h}
\end{eqnarray}
provided  $C>1$ (which we may assume).  Thus $h$ is large on $[r_j-C, r_j+C]$, and our  kernel will correlate well with quasimodes in $S_C(r_j)$.

On the other hand, since $\tilde{h}$ decays rapidly away from $\pm r_j$, and $\chi\in C^\infty$ implies that $\hat{\chi}$ decays rapidly, the convolution $h$ decays away from $\pm r_j$ as well; to be precise, the estimate $h(s) \lesssim r_j^{-1/2}\big||s|-r_j\big|^{-3}$ will suffice for our purposes.  Since $h$ is even, we can take $s\geq 0$ without loss of generality, and estimate first for $s>r_j$
\begin{eqnarray*}
h(s) & = & \int_{-\infty}^\infty \tilde{h}(r)\hat{\chi}(s-r) dr\\
& \leq & ||\tilde{h}||_\infty \int_{-\infty}^{\frac{s+r_j}{2}} \hat{\chi}(s-r)dr + ||\hat{\chi}||_\infty \int_{\frac{s+r_j}{2}}^\infty \tilde{h}(r)dr\\
& \leq & r_j^{-1/2} \int_{u=\frac{s-r_j}{2}}^{\infty} \hat{\chi}(u)du +  ||\hat{\chi}||_\infty \int_{u=\frac{s-r_j}{2}}^\infty \tilde{h}(r_j+u)du \\
& \lesssim & r_j^{-1/2}|s-r_j|^{-3} + r_j^{-1/2}\exp(-|s-r_j|/2C)\\
& \lesssim & r_j^{-1/2} |s-r_j|^{-3}
\end{eqnarray*}
since $\hat{\chi}(u) \lesssim |u|^{-4}$ for some uniform constant.
Similarly if $0<s<r_j$ we have
\begin{eqnarray*}
h(s) & \leq & ||\tilde{h}||_\infty\int_{-\infty}^{-r_j/2}\hat{\chi}(s-r)dr + ||\hat{\chi}||_\infty\int_{-r_j/2}^{\frac{r_j+s}{2}} \tilde{h}(r)dr + ||\tilde{h}||_\infty \int_{\frac{r_j+s}{2}}^\infty \hat{\chi}(s-r) dr\\
& \lesssim & r_j^{-1/2}|s-r_j|^{-3}
\end{eqnarray*}

Therefore, we can estimate various spectral integrals that will be needed in the argument:  we clearly have
\begin{equation}\label{h}
\int_0^\infty h(s) ds \lesssim r_j^{-1/2}
\end{equation}
and moreover
\begin{eqnarray}
\int_0^\infty sh(s) ds & \lesssim & r_j^{-1/2}\left(\int_0^{r_j-1} \frac{s}{(r_j-s)^3} ds +  \int_{r_j-1}^{r_j+1} s ds + \int_{r_j+1}^\infty \frac{s}{(s-r_j)^3}ds\right)\nonumber\\
& \lesssim & r_j^{-1/2} \left( \int_{u=1}^{r_j} \frac{r_j-u}{u^3}du +r_j + \int_{u=1}^\infty \frac{r_j+u}{u^3} \right)du\nonumber\\
& \lesssim & r_j^{1/2}\label{sh(s)}
\end{eqnarray}
and similarly 
\begin{eqnarray}
\int_0^\infty s h(s)^2 ds & \lesssim & r_j^{-1}\left(\int_0^{r_j-1} \frac{s}{(r_j-s)^6} ds +  \int_{r_j-1}^{r_j+1} s ds + \int_{r_j+1}^\infty \frac{s}{(s-r_j)^6}ds\right)\nonumber\\
& \lesssim  &  1\label{sh(s)^2}
\end{eqnarray}

We set $k$ to be the radial kernel corresponding to $h$, which by (\ref{Selberg/HC}) is supported in the ball of radius $\tau$ in $\mathbb{H}$.  We can estimate $||k||_{L^2(\mathbb{H})}$ by unitarity of the Helgason Fourier transform, giving
\begin{eqnarray*}
||k||_{L^2}^2 & = & \int_0^\infty h(s)^2s\tanh(\pi s)ds\\
& \leq & \int_0^\infty h(s)^2 s ds \lesssim 1
\end{eqnarray*}
by (\ref{sh(s)^2}).
Since $k$ is compactly supported inside the ball of radius $\tau$, which is less than the radius of injectivity of $M$, we can periodicize to obtain
$$k_p(z) = \sum_{\gamma\in\Gamma} k(p^{-1} \gamma z)$$
as a function on $M$, and if $\mathcal{F}_\Gamma\subset \mathbb{H}$ is a fundamental domain for $M$, we have $k_p(z) = k(p^{-1}z)$ whenever $p,z\in \mathcal{F}_\Gamma$.

\vspace{.2in}

We wish to show that the microlocal lift $|\kappa_p|^2dVol$ of $|k_p|^2dz$ is concentrated on cotangent vectors pointing radially towards and away from the base point $p$.  Precisely, let $\mathcal{B}_{r_j}\subset S^*\mathbb{H} = KAN$ be given by
$$	\mathcal{B}_{r_j} = \{ 	KAn_u : |u|\leq Nc^{-1}r_j^{-1/2}	\}	\cap B(0, \tau)K$$
for a constant $N$ to be chosen later; here $B(0,\tau)K = S^*\{d(z,0)<\tau\}\subset S^*\mathbb{H}$, where $\tau$ was chosen above to be less than  the radius of injectivity of $M$.  
Thus, $\mathcal{B}_{r_j}$ is an $O(r_j^{-1/2})$-neighborhood of the union of geodesic segments through $0$ up to distance $\tau$, with the implied constant depending on $N$ and $c$.
Define \`a la Lemma~\ref{asymp to op}
$$\kappa := b(r_j)r_j^{-1/4} \sum_{|n|\leq L=\lfloor \sqrt{r_j}\rfloor} \frac{L -|n|}{L} R^n k$$
where $R: \phi^{(r)}_{2n}\mapsto \phi^{(r)}_{2n+2}$ and $R^{-1}:\phi^{(r)}_{2n}\mapsto \phi^{(r)}_{2n-2}$ are the unitary raising and lowering operators as in section~\ref{ml lifts}.  The prefactor $b(r_j) \sim \sqrt{3/2}$ normalizes $||\kappa||_{L^2(S^*\mathbb{H})} = ||k||_{L^2(\mathbb{H})}$.  
Note  also that, since the summation over $\Gamma$ acting on the left of $g$ commutes with the left-invariant  operators $R^n$, we have
\begin{eqnarray*}
\kappa_p (g) 
& = & b(r_j)r_j^{-1/4} \sum_{|n|\leq L=\lfloor \sqrt{r_j}\rfloor } \frac{L -|n|}{L} R^n k_p\\
& = & \sum_{\gamma\in\Gamma} \kappa(p^{-1} \gamma.g)
\end{eqnarray*}
Additionally, since the unitary, left-invariant operators $R^n$ descend to $M$, and the weight spaces are orthogonal, we have 
$$||\kappa_p||_{L^2(S^*M)} = ||k_p||_{L^2(M)} = ||k||_{L^2(\mathbb{H})} \lesssim 1$$

Since the distribution
$$\sum_{n=-\infty}^\infty R^nk (g) = \int_{s=0}^\infty e^{(is-\frac{1}{2})\varphi(g.k_\theta)}d\theta h(s) s \tanh{(\pi s)} ds$$
by \cite{Z2},
the function $\kappa\in C^\infty(SL(2,\mathbb{R}))$ is given by
\begin{eqnarray*}
\kappa(g) & = & b(r_j)\int_{s=0}^\infty \frac{1}{2\pi}\int_{\theta=0}^{2\pi} e^{(is-\frac{1}{2})\varphi(g.k_\theta)}r_j^{-1/4}F_{\lfloor \sqrt{r_j}\rfloor}(\theta) d\theta h(s) s\tanh{(\pi s)} ds
\end{eqnarray*}
where\footnote{The standard normalization of the Fej\'er kernel is $||F_L||_{L^1}=1$, whereas our vectors are $L^2$-normalized, which causes a number of $r^{1/4}$ factors to appear throughout the discussion.} $F_{L}= \frac{1}{L}\left(\frac{\sin {L\theta/2}}{\sin{\theta/2}}\right)^2$ is the Fej\'er kernel of order $L$, and $\varphi\begin{pmatrix} a & b\\ c& d\end{pmatrix} = \log(a^2 + c^2)$ as in section~\ref{harmonic analysis}.

\vspace{.2in}

The next two Lemmas establish the key property, that $\kappa_p$ mainly lives in $p\mathcal{B}_{r_j}$.

\begin{Lemma}\label{phase derivative}
Let $g  =k_\alpha \begin{pmatrix} e^{t/2} & 0\\ 0 & e^{-t/2} 	\end{pmatrix}\begin{pmatrix} 1 & n \\ 0 & 1\end{pmatrix}\in B(0,\tau)K \backslash \mathcal{B}_{r_j}$; in particular,  $|t|\leq \tau$ and $Nc^{-1}r_j^{-1/2}\leq |n| \lesssim_\tau 1$.  Then 
$$\frac{d}{d\theta}\varphi(g.k_\theta) \gtrsim n$$
for all $|\theta|\leq |n|N^{-1/4}$.
\end{Lemma}

{\em Proof:} Since $\varphi$ is left $K$-invariant, it is sufficient to prove this for $g=\begin{pmatrix} e^{t/2} & 0\\ 0 & e^{-t/2} 	\end{pmatrix}\begin{pmatrix} 1& n\\ 0& 1\end{pmatrix} = \begin{pmatrix} e^{t/2} & e^{t/2}n\\ 0 & e^{-t/2} 	\end{pmatrix}$.  Write
\begin{eqnarray*}
g.k_\theta & =  & \begin{pmatrix}  e^{t/2} & e^{t/2}n\\ 0 & e^{-t/2} 	\end{pmatrix} \begin{pmatrix}  \cos{\theta} & -\sin{\theta}\\ \sin{\theta} & \cos{\theta} \end{pmatrix}\\
& = & \begin{pmatrix} e^{t/2} \cos{\theta} + e^{t/2} n \sin{\theta} & -e^{t/2}\sin{\theta}+e^{t/2}n\cos{\theta}\\ e^{-t/2}\sin{\theta} & e^{-t/2}\cos{\theta} \end{pmatrix}
\end{eqnarray*}
so that we wish to evaluate
\begin{eqnarray*}
\frac{d}{d\theta}\varphi(g.k_\theta)& = & \frac{d}{d\theta} \log\big( (e^{t/2}\cos{\theta} + e^{t/2}\sin{\theta})^2 + (e^{-t/2}\sin{\theta})^2\big)\\
& = & \frac{d}{d\theta}\log(e^t\cos^2{\theta} + ne^t\sin{2\theta} +e^tn^2\sin^2{\theta}+e^{-t}\sin^2{\theta})\\
& = & \frac{(-4\sinh{t}+e^tn^2)\sin{2\theta} + 2ne^{t}\cos{2\theta}}{e^t\cos^2{\theta} + O_{n,t}(\theta)}
\end{eqnarray*}
Since $|\theta|\leq |n|N^{-1/4}$, and $n$ and $t$ are uniformly bounded (depending only on $\tau$), we have $|(-4\sinh{t}+e^tn^2)\sin{2\theta}| \lesssim_\tau |n|N^{-1/4}$.  Since moreover the denominator is bounded below by $\frac{1}{2}e^{-\tau}$,  
the Lemma holds as soon as $N$ is large enough relative to $\tau$.  $\Box$

\begin{Lemma}\label{microlocalized kernel}
With notations as above, we can choose $N$ sufficiently large (independent of $r_j$) so that $\mathcal{B}_{r_j}$ satisfies
$$\int_{X\backslash p\mathcal{B}_{r_j}} |{\kappa}_p(y)|^2 dy \leq \frac{1}{50,000} c$$
\end{Lemma}

{\em Remark:} Recall that
$$||\kappa_p||_{L^2(S^*M)}^2 = ||k||_{L^2(\mathbb{H})}^2 \asymp 1$$
So Lemma~\ref{microlocalized kernel} is essentially saying that by letting $N$ be sufficiently large (independent of $r_j$), we can get a large percentage of the $L^2$-mass of $\kappa_p$ to lie inside $p\mathcal{B}_{r_j}$.

{\em Proof:}  
The main step is bounding $\int |\kappa|^2$ on $B(0,\tau)K\backslash\mathcal{B}_{r_j}$.  For this, we will decompose $\kappa$ into two pieces--- the contribution of low frequencies, and that of high frequencies--- and estimate each individually on $B(0,\tau)K\backslash\mathcal{B}_{r_j}$.

Write $\kappa = \kappa_1 + \kappa_2$, where 
\begin{eqnarray*}
\kappa_1 (g) & = & b(r_j) \int_{s=0}^{\eta r_j}\frac{1}{2\pi}\int_{\theta=0}^{2\pi} e^{(is-\frac{1}{2})\varphi(g.k_\theta)}r_j^{-1/4}F_{\lfloor\sqrt{r_j}\rfloor}(\theta) d\theta h(s) s\tanh{(\pi s)} ds\\
\kappa_2 (g) & = & b(r_j) \int_{s=\eta r_j}^\infty \frac{1}{2\pi}\int_{\theta=0}^{2\pi} e^{(is-\frac{1}{2})\varphi(g.k_\theta)}r_j^{-1/4}F_{\lfloor\sqrt{r_j}\rfloor}(\theta) d\theta h(s) s\tanh{(\pi s)} ds
\end{eqnarray*}
where $\eta$ is chosen small enough, depending on $C$ and $c$, so that computing as in (\ref{sh(s)^2}) 
\begin{eqnarray*}
 \int_{s=0}^{\eta r_j} h(s)^2 s \tanh{(\pi s)}ds & \leq & r_j^{-1}\int_{s=0}^{\eta r_j} \frac{s}{(r_j-s)^6}ds\\
 & \leq & r_j^{-1} \int_{(1-\eta)r_j}^{r_j} \frac{r_j-u}{u^6} du\\
 & < & \frac{1}{200,000} c
 \end{eqnarray*}
Thus 
\begin{eqnarray}
\int_{B(0,\tau)K\backslash\mathcal{B}_{r_j}} |\kappa_1|^2 & \leq & ||\kappa_1||_{L^2(S^*\mathbb{H})}^2\leq ||k_1||^2_{L^2(\mathbb{H})} \nonumber\\
& \leq &   \int_{s=0}^{\eta r_j} h(s)^2 s \tanh{(\pi s)}ds < \frac{1}{200,000} c \label{kappa1}
\end{eqnarray}
for $k_1(z)= \int_{s=0}^{\eta r_j} \frac{1}{2\pi}\int_{\theta=0}^{2\pi} e^{(is-\frac{1}{2})\varphi(z.k_\theta)}d\theta h(s) s\tanh{(\pi s)} ds$, by unitarity of the Helgason Fourier transform.  

We turn to  $\kappa_2$ on $B(0,\tau)K\backslash\mathcal{B}_{r_j}$; we will estimate $|\kappa_2|$ pointwise.  For any $s > \eta r_j$, and any $g\in KA\begin{pmatrix}1&n\\0&1\end{pmatrix}$ with $|n|\geq Nc^{-1}r_j^{-1/2}$,  write
\begin{eqnarray}
\lefteqn{\left|\int_{\theta=0}^{2\pi} e^{(is-\frac{1}{2})\varphi(g.k_\theta)}F_{\lfloor\sqrt{r_j}\rfloor}(\theta) d\theta \right|}\nonumber\\
& \leq & \int_{|\theta|\geq |n|/N^{1/4}} \left|F_{\lfloor\sqrt{r_j}\rfloor}(\theta)  \right| e^{-\frac{1}{2}\varphi(g.k_\theta)}d\theta
+ \sum_{|m|\leq \sqrt{r_j}}\left|\int_{|\theta|\leq |n|/{N}^{1/4}}  e^{i(s \varphi(g.k_\theta) +m\theta)} e^{-\frac{1}{2}\varphi(g.k_\theta)} d\theta\right|\nonumber\\
& & \label{pre kappa2}
\end{eqnarray}
Now since $e^{-\frac{1}{2}\varphi(g.k_\theta)}$ and its derivatives are uniformly bounded on $B(0,\tau)K$, we may apply a non-stationary phase argument to determine that each integral in the sum on the right is $O(s^{-1}|n|^{-1})$ whenever the derivative $s\varphi'(g.k_\theta) + m$ of the phase function is $ \gtrsim s |n|$: setting $a(\theta) = e^{-\frac{1}{2}\varphi(g.k_\theta)}$, and recalling that $\varphi$ and $a$, and their derivatives,  are uniformly bounded on $B(0,\tau)K$, we write
\begin{eqnarray*}
\lefteqn{\left|\int_{|\theta|\leq |n|/{N}^{1/4}}  e^{i(s \varphi(g.k_\theta) +m\theta)} a(\theta) d\theta\right|}\\ 
& = & \left|\int_{|\theta|\leq |n|/{N}^{1/4}}  \frac{1}{i(s\varphi'(g.k_\theta)+m)} \frac{d}{d\theta}\left[e^{i(s \varphi(g.k_\theta) +m\theta)}\right] a(\theta) d\theta\right|\\
& = & \left|\int_{|\theta|\leq |n|/{N}^{1/4}} \frac{d}{d\theta}\left[e^{i(s \varphi(g.k_\theta) +m\theta)}\right] \frac{a}{i(s\varphi'(g.k_\theta)+m)}d\theta\right|\\
& \leq & 2\sup \left|\frac{a}{s\varphi'(g.k_\theta)+m} \right| 
+ \int_{|\theta|\leq |n|/N^{1/4}} \left| \frac{a'(s\varphi' +m) - a s\varphi''}{(s\varphi'+m)^2}\right| d\theta\\
& \lesssim & s^{-1}|n|^{-1} + \int_{|\theta|\leq |n|/N^{1/4}} \frac{s}{(s\varphi'+m)^2} d\theta\\
& \lesssim & s^{-1}|n|^{-1} +  s^{-1}|n|^{-2} \int_{|\theta|\leq |n|/N^{1/4}} d\theta\\
& \lesssim &  s^{-1}|n|^{-1} 
\end{eqnarray*}
But Lemma~\ref{phase derivative} shows that $s \frac{d\varphi(g.k_\theta)}{d\theta} \gtrsim s |n| > \eta Nc^{-1}r_j^{1/2} > 2|m|$ for all $|\theta|\leq |n|/N^{1/4}$ if $N$ is chosen large enough, so  that indeed $s\varphi' + m \gtrsim s|n|$.  Therefore
\begin{eqnarray*}
\sum_{|m|\leq \sqrt{r_j}} \left|\int_{|\theta|\leq |n|/{N}^{1/4}}  e^{i(s \varphi(g.k_\theta) +m\theta)} e^{-\frac{1}{2}\varphi(g.k_\theta)} d\theta\right| & \lesssim & \sum_{|m|\leq\sqrt{r_j}}  s^{-1}|n|^{-1}\\
& \lesssim & r_j^{1/2}s^{-1}|n|^{-1}
\end{eqnarray*} 

To estimate the $|F_{\lfloor \sqrt{r_j}\rfloor}|$ term of (\ref{pre kappa2}), observe that since $\varphi$ is bounded on $B(0,\tau)K$
\begin{eqnarray*}
\int_{|\theta|\geq |n|/N^{1/4}} |F_{\lfloor\sqrt{r_j}\rfloor}(\theta)|e^{-\frac{1}{2}\varphi(g.k_\theta)}d\theta & \lesssim & \int_{|\theta|\geq |n|/N^{1/4}} |F_{\lfloor\sqrt{r_j}\rfloor}(\theta)|d\theta\\
& \lesssim  & \int_{|\theta|\geq |n|/{N}^{1/4}} \frac{1}{\sqrt{r_j}} \left(\frac{\sin{(\lfloor\sqrt{r_j}\rfloor\theta/2)}}{\sin{\frac{\theta}{2}}}\right)^2 d\theta\\
& \lesssim & \int_{|\theta|\geq |n|/N^{1/4}} r_j^{-1/2} \theta^{-2} d\theta \lesssim {N}^{1/4} r_j^{-1/2}|n|^{-1}
\end{eqnarray*}
so that combining the two parts of (\ref{pre kappa2}) we have
\begin{equation}
\left|\int_{\theta=0}^{2\pi} e^{(is-\frac{1}{2})\varphi(g.k_\theta)}F_{\lfloor\sqrt{r_j}\rfloor}(\theta) d\theta \right| \leq N^{1/4}|n|^{-1} (r_j^{-1/2} + r_j^{1/2}s^{-1})
\end{equation}

Therefore, since 
\begin{eqnarray*}
|\kappa_2(g)| & \leq & \int_{s=\eta r_j}^\infty \frac{1}{2\pi}\left|\int_{\theta=0}^{2\pi} e^{(is-\frac{1}{2})\varphi(g.k_\theta)}F_{\lfloor\sqrt{r_j}\rfloor}(\theta) d\theta\right| r_j^{-1/4} h(s) s\tanh{(\pi s)} ds\\
& \lesssim & N^{1/4}|n|^{-1}r_j^{-1/2} \int_{s=\eta r_j}^\infty (r_j^{-1/4} sh(s) + r_j^{3/4}h(s)) ds  \\
& \lesssim & N^{1/4} r_j^{-1/4} |n|^{-1}
\end{eqnarray*}
by (\ref{h}) and (\ref{sh(s)}),  we see that
\begin{eqnarray*}
\int_{B(0,\tau)K\backslash\mathcal{B}_{r_j}} |\kappa_2|^2 & \leq & \int_{|n|\geq Nc^{-1}r_j^{-1/2}} |\kappa_2|^2\\
& \lesssim & N^{1/2} r_j^{-1/2} \int_{|n|\geq Nc^{-1}r_j^{-1/2}} |n|^{-2}\\
& \lesssim & N^{1/2} r_j^{-1/2} \cdot (N^{-1} c r_j^{1/2})\\
& \lesssim & N^{-1/2} c
\end{eqnarray*}  
and we now choose $N$ large enough (depending on the implied constant, which in turn depends on the parameters $\tau$, $C$, and $c$) so that
\begin{equation}
\int_{B(0,\tau)K\backslash\mathcal{B}_{r_j}} |\kappa_2|^2 < \frac{1}{200,000}c  \label{kappa2}
\end{equation}
Combining (\ref{kappa1}) and (\ref{kappa2}) we get
\begin{eqnarray*}
\int_{S^*M\backslash p\mathcal{B}_{r_j}} |\kappa(p^{-1} g)|^2 dg & = & \int_{B(0,\tau)K\backslash \mathcal{B}_{r_j}} |\kappa(g)|^2 dg\\
& \leq & 2 \int_{B(0,\tau)K\backslash\mathcal{B}_{r_j}}(|\kappa_1|^2 + |\kappa_2|^2)dg\\
& \leq & \frac{1}{50,000}c
\end{eqnarray*}
as required.  $\Box$

\section{Proof of Theorem~\ref{optimality}}\label{optimal quasimodes}
Equipped with this construction, we can return to the main result.

{\em Proof of Theorem~\ref{optimality}:}
The next order of business is to find a quasimode $\psi_j \in \mathcal{S}_j$, and a point $p_j \in \Gamma\backslash \mathbb{H}$ at which $|{\psi_j}|$ is large.  The following observation can be found in \cite{SarnakLetterMorawetz}, in the context of constructing an eigenfunction of large $L^\infty$-norm in a highly degenerate eigenspace.

Consider 
$$\int_{M} \sum_{\phi_l\in \mathcal{S}_j} |\phi_l(z)|^2 dx = \sum_{\phi_l\in \mathcal{S}_j} \int_{M} |\phi_l(z)|^2 dx = \sum_{\phi_l \in \mathcal{S}_j} 1 =c r_j$$
This implies that there exists a point $p_j\in M$, such that 
$$\sum_{\phi_l\in \mathcal{S}_j} |\phi_l(p_j)|^2 \geq c r_j $$
is at least as large as the average value.  Therefore the quasimode
$$\psi_j := \sum_{\phi_l \in \mathcal{S}_j} {\phi_l(p_j)}\phi_l$$
satisfies
\begin{equation}\label{massive point}
\psi_j(p_j) = ||\psi_j||_2^2 = \sum_{\phi_l \in \mathcal{S}_j} |{\phi_l(p_j)}|^2   \geq \sqrt{cr_j}||\psi_j||_2
\end{equation}
This is the sequence $\psi_j$ of quasimodes we will use.

\vspace{.2in}

Now consider $\mathcal{B}_{j}:=  p_j\mathcal{B}_{r_j}$, and since 
$$||{\kappa}_{p_j}||^2_{L^2(\mathcal{B}_{j})} \leq ||\kappa_{p_j}||^2_{L^2(S^*M)}= ||k_{p_j}||^2_{L^2(M)} \lesssim 1$$ 
we have by Cauchy-Schwarz
\begin{eqnarray*}
\int_{\mathcal{B}_{j}} |{\Psi_j}(x)|^2 dx & \geq   & \frac{1}{||\kappa_{p_j}||_{L^2(\mathcal{B}_j)}^2} \left(\int_{\mathcal{B}_{j}} |{\Psi_j}(x)| |{\kappa}_{p_j} (x)| dx\right)^2 \\
& \gtrsim &  \Big( \langle {\Psi_j}, {\kappa}_{p_j}\rangle - ||{\kappa}_{p_j}||_{L^2(X\backslash\mathcal{B}_j)}||\Psi_j||_2 \Big)^2\\
& \gtrsim & \Big( \langle {\Psi_j}, {\kappa}_{p_j}\rangle - \frac{1}{200} \sqrt{c}||\Psi_j||_2 \Big)^2
\end{eqnarray*}
since $\int_{X\backslash \mathcal{B}_j} |{\kappa}_{p_j}|^2 \leq \frac{1}{50,000} c $ by Lemma~\ref{microlocalized kernel}.

Now since the normalized raising and lowering operators $R^n$ are unitary, and the weight spaces orthogonal, we have $\langle \Psi_j, \kappa_{p_j}\rangle = \langle \psi_j, k_{p_j}\rangle$, so that writing the spectral expansion $k_{p_j}(z) = \sum_{\phi_l} h(s_l)\phi_l(p_j)\phi_l(z)$, we use (\ref{big h}) to get
\begin{eqnarray*}
\langle \Psi_j, \kappa_{p_j} \rangle & = & \langle \psi_j, k_{p_j} \rangle \\
& = & \sum_{\phi_l\in \mathcal{S}_j} h(s_l) |\phi_l(p_j)|^2\\
& \geq &  \frac{1}{100}r_j^{-1/2} \sum_{\phi_l\in \mathcal{S}_j} |\phi_l(p_j)|^2\\
& \geq &  \frac{1}{100}r_j^{-1/2} \cdot (cr_j)^{1/2} ||\Psi_j||_2 = \frac{1}{100} \sqrt{c}||\Psi_j||_2
\end{eqnarray*}
Therefore
\begin{eqnarray*}
\int_{\mathcal{B}_{j}} |{\Psi_j}(x)|^2 dx & \gtrsim  &  \left(\left[\frac{1}{100}-\frac{1}{200} \right]\sqrt{c} ||{\Psi_j}||_2\right)^2 \gtrsim  ||{\Psi_j}||_2^2
\end{eqnarray*}
with implied constant depending on $c$, $C$, and $\tau$.  In other words, the $L^2$-mass of ${\Psi_j}$ inside $\mathcal{B}_{j}$ is at least a fixed percentage of the total mass, independent of $r_j\to\infty$.

To finish, suppose we have a subsequence of $\{\Psi_j\}$ such that $d\mu_j = |\Psi_j|^2dVol$ converges weak-* to a measure $\mu$, and pick a further subsequence of the $j$'s so that $p_j\to \bar{p}\in \Gamma\backslash\mathbb{H}$.  Consider any fixed neighborhood $U\subset X$ of the (compact) union of geodesic segments of length $\tau$ through $\bar{p}$; it is evident that $U$ can be chosen to have arbitrarily small volume in $X$.  On the other hand, any such neighborhood must contain $\mathcal{B}_{j}$ for sufficiently large $j$ in our subsequence, whereby $\mu_j(U) = \int_U |\Psi_j(x)|^2 dx \gtrsim ||\Psi_j||_2^2$ for all sufficiently large $r_j$.  Thus the measure $\mu$ concentrates a positive proportion of its mass on this codimension $1$ subset.   $\Box$

\def\cprime{$'$}

\end{document}